\documentclass{amsart}

\usepackage{pstricks, pst-tree, pst-node} 
\usepackage[latin1]{inputenc}
\usepackage{amssymb}
\usepackage{xfrac}     
\usepackage{lmodern}
\usepackage{stmaryrd}
\usepackage{graphicx}
\usepackage{amsmath}  
\usepackage{amsthm}
\usepackage{wrapfig} 
\usepackage{color}  
\usepackage{microtype} 
\usepackage{setspace}
\usepackage{amsfonts}
\usepackage{wasysym}
\usepackage{cancel}
\usepackage{verbatim}
\usepackage{dsfont}
\usepackage{paralist}

  \usepackage{nicefrac}

\usepackage{diagbox}

\newcommand{\F}{\mathbb{F}}
\newcommand{\Pcal}{\mathcal{P}}
\usepackage{enumerate, verbatim, cleveref}

 \DeclareMathOperator{\symlen}{\mathrm{sl}}
\newcommand{\Gq}{\mathsf{G}}
\newcommand{\D}{\mathsf{D}}
\renewcommand{\P}{\mathsf{P}}
\newcommand{\GP}{\mathsf{GP}}
\newcommand{\W}{\mathsf{W}}

\DeclareMathOperator{\I}{I}

\DeclareMathOperator{\Z}{\mathbb{Z}}
\DeclareMathOperator{\Q}{\mathbb{Q}}

\DeclareMathOperator{\N}{\mathbb{N}}

\DeclareMathOperator{\CC}{\mathbb{C}}

\newcommand{\spn}{\mathrm{span}}

\newcommand{\qf}[1]{\langle #1\rangle}

\newcommand{\laurent}[1]{(\!(#1)\!)}
\newcommand{\Pfister}[1]{\langle\!\langle #1\rangle\!\rangle}

\newtheoremstyle{plain2}
  {10pt}   
  {10pt}   
  {\itshape}  
  {0pt}       
  {\bfseries} 
  {}         
  {5pt plus 1pt minus 1pt} 
  {}          
  
 \newtheoremstyle{beweis}
  {10pt}   
  {10pt}   
  {\normalfont}  
  {0pt}       
  {\bfseries} 
  {:}         
  {5pt plus 1pt minus 1pt} 
  {}          

\newtheoremstyle{definition2}
  {10pt}   
  {10pt}   
  {\normalfont}  
  {0pt}       
  {\bfseries} 
  {}         
  {5pt plus 1pt minus 1pt} 
  {}          

\theoremstyle{plain2}
\newtheorem{satz}{Satz}[section]
\newtheorem{lemma}[satz]{Lemma}
\newtheorem{proposition}[satz]{Proposition}
\newtheorem{theorem}[satz]{Theorem}
\newtheorem{corollary}[satz]{Corollary}

\theoremstyle{definition2}
\newtheorem{definition}[satz]{Definition}

\newtheorem{remark}[satz]{Remark}
\newtheorem{example}[satz]{Example}

\theoremstyle{beweis}

\begin{document}

\title{On the Number of Pfister Forms and the Symbol Length of Fields With Finite Square Class Number}

\author{Detlev Hoffmann}
\address{Fakult\"at f\"ur Mathematik, Technische Universit\"at Dortmund, D-44221 Dortmund, Germany}
\email{detlev.hoffmann@tu-dortmund.de}

\author{Nico Lorenz}
\address{Fakult\"at f\"ur Mathematik, Ruhr-Universit\"at Bochum, Universit\"atsstra\ss e 150, 44801 Bochum, Deutschland}
\email{nico.lorenz@ruhr-uni-bochum.de}
\date{\today}



\begin{abstract}
Let $F$ be a field of characteristic not $2$ with finitely many square classes. 
Using  combinatorial arguments applied to objects related to vector spaces over finite fields, we deduce an upper bound for the number of Pfister forms over $F$. 
Moreover, we compute upper bounds for the $n$-symbol length of $F$ ($n\in\N$),
i.e., the smallest
integer $\symlen_n(F)\geq 0$ such that to each quadratic form
$\phi\in \I^n(F)$ there exists some $0\leq k\leq \symlen_n(F)$ and
Pfister forms $\pi_1,\ldots, \pi_k$ such that
$\varphi\equiv \pi_1+\ldots+\pi_k\mod \I^{n+1}(F)$.
In particular, we rediscover a bound that can also be deduced from
a result by Bruno Kahn that he stated without proof.

Classification (MSC 2020): 11E04, 11E81, 12D15, 12G05, 19D45

\end{abstract}

\maketitle

\section{Introduction}

Throughout this paper, let $F$ be a field of characteristic different from $2$. 
By a \emph{quadratic form} or just \emph{form} for short, we will
always mean a finite dimensional non-degenerate quadratic form over $F$.

An $n$\emph{-fold Pfister form} $\pi$
for some $n\in\N$ is an $n$-fold tensor product of
binary forms that represent $1$, i.e., $\pi$ is a  form of the shape
$\Pfister{a_1,\ldots, a_n}:=\qf{1, -a_1}\otimes\ldots\otimes\qf{1, -a_n}$
with $a_1,\ldots, a_n\in F^\ast=F\setminus\{ 0\}$. 
The set of $n$-fold Pfister forms over $F$ is denoted by $\P_n(F)$;
the set of forms that are similar to $n$-fold Pfister forms is
denoted by $\GP_n(F)$. 
Forms in $\GP_n(F)$ are called \emph{general Pfister forms}.

Both $\P_n(F)$ and $\GP_n(F)$ generate additively the $n$-th power $\I^n(F)$ of the
\emph{fundamental ideal} $\I(F)$ of Witt classes of even-dimensional forms
in the Witt ring $\W(F)$, i.e., any Witt class $\varphi\in \I^n(F)$ can be
written as a sum $\varphi=\pi_1+\ldots+\pi_m$ for some $m$ and with
$\pi_i\in \GP_n(F)$ (resp. with $\pi_i$ or $-\pi_i\in \P_nF$).
By convention, we denote by $\qf1$ the unique $0$-fold Pfister form and
set $\I^0(F)=\W(F)$.

For a form $\varphi\in \I^n(F)$, we are interested in the
\emph{$n$-symbol length} of $\varphi$ (or \emph{symbol length} for short
if the integer is clear from the context). The $n$-symbol length is defined as
\[\symlen_n(\varphi):=\min\{k\in\N\mid \exists \pi_1,\ldots, \pi_k\in \GP_n(F): 
\varphi\equiv \pi_1+\ldots+\pi_k\mod \I^{n+1}(F)\}.\]
Note that one may have used $\P_n(F)$ in this definition instead since for all
$c\in F^*$ and all $\pi\in\P_n(F)$, one has $\pi\equiv c\pi\mod \I^{n+1}(F)$.

We further define 
\[\symlen_n(F):=\sup_{\varphi\in \I^n(F)}\{\symlen_n(\varphi)\}\in\N\cup\{\infty\}.\]
(In the literature, one also finds the notation
$\lambda_n(F)$ or $\lambda^n(F)$ to denote the symbol length.) 

Due to the groundbreaking results by V.~Voevodsky et al.
in \cite{MR2276765} and \cite{MR2031199}, the symbol length has
connections to Milnor $K$-theory and Galois cohomology.  In characteristic not $2$, 
the Milnor $K$-groups $K_n(F)/2$, the Galois cohomology groups  $H^n(F,\Z/2\Z)$
and $\I^n(F)/\I^{n+1}(F)$ are all mutually isomorphic.  Under these isomorphisms,
generators of the form $\{a_1,\ldots, a_n\}$ in $K_n(F)/2$ correspond
to $n$-fold cup products $(a_1)\cup\ldots\cup(a_n)$ in $H^n(F,\Z/2\Z)$ and
to $n$-fold Pfister forms $\Pfister{a_1,\ldots,a_n}$ modulo $\I^{n+1}(F)$, respectively.
So the analogous notions of symbol lengths can be translated from any one of these
three groups to any other.

The determination of the symbol length is an active field of research,
see, e.g., \cite{MR3404396}, \cite{MR3413868} or \cite{MR4102695}
for some recent progress in the field.  The $3$-symbol length (and the related
concept of $3$-Pfister number, i.e., the minimal length of a representation of
forms in  $I^3(F)$ as sums of general $3$-fold Pfistter forms)  has, for example,
been studied in \cite{MR2630047}, \cite{MR2537109}, \cite{MR3010541}, \cite{MR4554566}.

The aim of this article is to find upper bounds 
for $\symlen_n(\varphi)$ and $\symlen_n(F)$ in the case where the square
class group $F^*/F^{*2}$ is finite, i.e., $q(F)=|F^*/F^{*2}|=2^d<\infty$.
As this problem is rather trivial for  $n=1$ and in order to avoid
subtleties in our calculations, we assume throughout that $n$ be an
integer $\geq 2$.
Note that K.~Becher and the first author obtained best possible general upper bounds
in the case $n=2$ in terms of the cardinality $q(F)$, see
\cite{MR2061565}.

Our estimates depend on some field invariants such as the level,
the Pythagoras number and the size of certain quotients of
subgroups of the multiplicative group
of the field related to sums of squares. 

In Section 2, we provide a counting argument for the number of
Pfister forms of a certain type over a field with finitely many square classes,
and provide upper bounds for the size of certain quotient groups
in a filtration of the square class group of a field with finitely many
square classes.  Using these results, we use more refined counting methods to
obtain bounds for the $n$-symbol length (\Cref{PolynomialBound}).
The bound is of a rather technical nature and can be weakened to yield a
polynomial bound in $d$ of degree $n-1$ (\Cref{kahn-bound}). Such a polynomial
can be computed explicitly (see \Cref{example-n=3} where such polynomials are
explicitly computed in the case $n=3$).
A polynomial bound of that type
can also be deduced from a result that was stated by B.~Kahn but
without proof, \cite[Proposition 2.3(h)]{Kahn}.

In the remainder of this introduction, we recall some further basic definitions
and facts from the algebraic theory of quadratic forms and refer to \cite{Lam2005}
for any undefined terminology or additional facts that we state without
further reference.

Isometry of two forms $\varphi_1,\varphi_2$ will be denoted by
$\varphi_1\cong\varphi_2$, their orthogonal sum by $\phi_1\perp\phi_2$, and their
tensor product by $\phi_1\otimes\phi_2$.  The orthogonal sum of $m$ copies of $\phi$
will be written $m\times\phi$.
By abuse of notation, we will denote the \emph{Witt class} of a quadratic form
$\varphi$ in the \emph{Witt ring} $\W(F)$ again by $\varphi$.

For a quadratic form $\varphi$ defined over a vector space $V$,
we denote by $\D_F(\varphi)=\{a\in F^\ast\mid \exists v\in V: \varphi(v)=a\}$
the set of nonzero elements represented by $\varphi$ and
by $\Gq_F(\varphi)=\{a\in F^\ast\mid a\varphi\cong\varphi\}$ the multiplicative
group of similarity factors of $\varphi$.
We define $\D_F(m)=\D_F(m\times\qf{1})$, the nonzero sums of $m$ squares in $F$.

A form $\phi$ is called \emph{round} if 
$\D_F(\phi)=\Gq_F(\phi)$.  Pfister forms are always round.  In particular,
if $\varphi$ is isometric to the $m$-fold Pfister form
$\Pfister{-1,\ldots, -1}$, we obtain $\Gq_F(\varphi)=\D_F(\varphi)=\D_F(2^m)$,
the group of nonzero sums of $2^m$ squares in $F$.
The set of nonzero sums of squares in $F$ will be denoted by
$\D_F(\infty)=\bigcup_{m=1}^{\infty}\D_F(m)$.

The \emph{level} of $F$, denoted by $s(F)$, is the least integer
$k$ such that $-1$ is a sum of $k$ squares in $F$, or $\infty$ if no such $k$ exists.
A famous result by A.~Pfister \cite{MR0175893} states that
$s(F)$ is a power of $2$ if finite, and that each power of $2$ occurs as
a level of a suitable field.
A field is called \emph{formally real} (or \emph{real} for short),
if $-1$ is not a sum of squares, and \emph{nonreal} otherwise.
The \emph{Pythagoras number} of $F$, denoted by $p(F)$ is the least
integer $k$ such that any sum of squares in $F$ is a sum of $k$ squares,
or $\infty$ if no such integer exists. 
For nonreal fields, we have $s(F)\leq p(F)\leq s(F)+1$,
see \cite[Chapter XI. Theorem 5.6]{Lam2005}.
For real fields, it was shown by the first author in \cite[Theorem 1]{MR1670858} that 
every positive integer can be realized as the Pythagoras number of a suitable field.

Since we will focus on fields with finite square class number,
it should be noted that there is no known example of a field with $q(F)<\infty$ and
finite level $s(F)\geq 8$.  In fact, if the so-called
\emph{Elementary Type Conjecture} were true, then any field with
$q(F)<\infty$ would have level $s(F)\in\{ 1,2,4,\infty\}$, see also 
\cite[Chapter XIII. Question 6.2]{Lam2005}.

We will use several combinatorial arguments in the sequel. We need
the \emph{binomial coefficient} $\binom nk$ for integers $n,k\in\Z$, i.e.,
the number of subsets with $k$ elements of a set with $n$ elements. 
For $0\leq k\leq n$, we have $\binom nk=\frac{n!}{k!(n-k)!}$ and in all other
cases, we have $\binom nk=0$.
Finally, for a set $S$, we denote the \emph{power set} of $S$, i.e.,
the set consisting of all subsets of $S$ by $\Pcal(S)$.
More generally, for $m\in\N$, we denote by $\Pcal_m(S)$ the set of all
subsets $A\subseteq S$ with $|A|=m$.
In particular, if $S$ is finite, then $|\Pcal_m(S)|=\binom{|S|}m$.

\section{The Number of Pfister Forms}\label{Sec3}

In this section, we obtain an upper bound for the number of isometry
classes of Pfister forms over fields with finite square class number.
Since we will use a crucial relation between Pfister forms and
certain vector spaces, we collect some basic facts from linear algebra. 

Let $\F_q$ be the finite field with $q$ elements, let $0\leq m\leq d<\infty$ and
denote by $V(q,d,m)$ the number of $m$-dimensional subspaces
of a $d$-dimensional vector space over $\F_q$
This number is well known (or an easy exercise in linear algebra)
and given in the following result.

\begin{proposition}\label{prop:NumberOfSubspaces}
	We have
	\[V(q,d,m)=\prod_{\ell=0}^{m-1}\frac{q^{d-\ell}-1}{q^{m-\ell}-1}.\]
\end{proposition}

We will apply these results in particular to the groups
$F^\ast/F^{\ast2}, F^\ast/\D_F(2^m)$ and $F^\ast/\pm\D_F(2^m)$
for some nonnegative integer $m$.
Since these groups all have exponent $\leq2$, they can be interpreted
as vector spaces over $\mathbb F_2$. 
Thus, if finite, these groups have order $2^k$ and thus
$\mathbb F_2$-dimension $k$ for some nonnegative integer $k$.

We therefore give an upper bound in the special case of
$\mathbb F_2$-vector spaces for later reference.

\begin{corollary}\label{cor:BoundNumberSubspaces}
$$V(2,d,m)\leq 2^{m(d-m)}\alpha_m\quad\mbox{where}\quad
\alpha_m=\prod_{\ell=1}^m\left(1+\frac{1}{2^\ell -1}\right).$$	
\end{corollary}
\begin{proof}
	Using \Cref{prop:NumberOfSubspaces} we have 
	\begin{align*}
		V(2,d,m)&=\prod_{\ell=0}^{m-1}\frac{2^{d-\ell}-1}{2^{m-\ell}-1}\\
		&=\prod_{\ell=0}^{m-1}\frac{2^{d-m}\left(2^{m-\ell}-1\right)+2^{d-m}-1}{2^{m-\ell}-1}\\
		&<\prod_{\ell=0}^{m-1}\frac{2^{d-m}\left(2^{m-\ell}-1\right)+2^{d-m}}{2^{m-\ell}-1}\\
		&=\prod_{\ell=0}^{m-1}2^{d-m}\cdot\left(1+\frac1{2^{m-\ell}-1}\right)\\
		&=2^{m(d-m)}\cdot\prod_{\ell=1}^{m}\left(1+\frac1{2^{\ell}-1}\right)\\
	\end{align*}
\end{proof}

\begin{remark} One readily obtains
$$\alpha_m=\left\{\begin{array}{rl}
1 & \mbox{if $m=0$;}\\
2 & \mbox{if $m=1$;}\\
\frac{8}{3} & \mbox{if $m=2$.}
\end{array}\right.$$
For $m\geq 3$, we write
$$\alpha_m=\frac{8}{3}\prod_{\ell=3}^m\left(1+\frac{1}{2^\ell -1}\right)$$
and note that for $\ell\geq 3$, we have $3\cdot 2^{\ell-2}\leq 2^\ell -1$, so
$$\alpha_m\leq\frac{8}{3}\prod_{\ell=1}^{m-2}\left(1+\frac{1}{3\cdot 2^\ell}\right)$$
Using $1+x\leq e^x$ we get for all $k\geq 1$,
$$\prod_{\ell=1}^k\left(1+\frac{1}{3\cdot 2^\ell}\right)\leq
\exp\left(\frac{1}{3}\sum_{\ell=1}^k2^{-\ell}\right)< e^{1/3}$$
and thus, $\alpha_m<\frac{8}{3}e^{1/3}$.
\end{remark}

Let $\P_{n,m}(F)$ be the set of isometry classes of \emph{anisotropic}
$n$-fold Pfister forms over $F$ that can be written with $m$ slots that
are $-1$, but such that there is no representation with $m+1$ slots that are $-1$.
Of course, this set may be empty. 

\begin{lemma}\label{LinearIndependentRep}
	Let $\varphi\in \P_{n,m}(F)$ and $x_1,\ldots, x_{n-m}\in F^\ast$ such that we have 
	\[\varphi\cong\Pfister{-1,\ldots, -1,-x_1,\ldots, -x_{n-m}}.\] 
	Then  $x_1, \ldots, x_{n-m}$ are linearly independent 
when considered as elements in the $\F_2$-vector space
${F^\ast/\pm \D_F(2^m)}$ and thus also in ${F^\ast/ \D_F(2^m)}$.
\end{lemma}
\begin{proof}
	If $x_1, \ldots, x_{n-m}$ are linearly dependent in ${F^\ast/\pm \D_F(2^m)}$, after renumbering, there is some $i\in\{0,\ldots, n-m-1\}$ and some $c\in\pm \D_F(2^m)$ with 
	\[x_{i+1}=c\prod_{k=1}^ix_k.\]
	We claim that we can find $a\in \pm \D_F(2^m)$ and
$n-m-1$ elements $y_1,\ldots, y_{n-m-1}$ in
$\{ x_1,\ldots,x_{n-m}\}$ such that we have 
	\[\varphi\cong\Pfister{-1,\ldots, -1, a, -y_1,\ldots, -y_{n-m-1}}\]
	This is clear for $i=0$, so let now $i\geq1$.
	Recall that for all $a,b\in F^\ast$, we have an isometry $\Pfister{-a,-b}\cong\Pfister{-a, -ab}$.
	Using this equality $i$ times, we obtain
	\begin{align*}
		\Pfister{-x_1,\ldots, -x_i, -x_{i+1}}&\cong
\Pfister{-x_1,\ldots, -x_i,-cx_1\cdot\ldots\cdot x_i}\\
		&\cong\Pfister{-x_1,\ldots, -x_i, -c}\cong
\Pfister{-c, -x_1,\ldots, -x_i},
	\end{align*}
	i.e. for $a=-c\in\pm \D_F(2^m)$ we have
    		\[\varphi\cong\Pfister{-1,\ldots, -1, a, -x_1,\ldots, -x_{i}, \widehat{-x_{i+1}}, -x_{i+2},\ldots, -x_{n-m-1}}\]
    	(here, the hat means to omit $-x_{i+1}$).  This establishes the claim.
	We now distinguish whether $a\in\D_F(2^m)$ or $a\in-\D_F(2^m)$.
	For $a\in \D_F(2^m)=\Gq_F(2^m)$ we have 
	\begin{align*}
       		\Pfister{-1,\ldots, -1, a}&=\Pfister{-1,\ldots, -1}\perp -a\Pfister{-1,\ldots, -1}\\
        		&\cong\Pfister{-1,\ldots, -1} \perp-\Pfister{-1,\ldots, -1},
    	\end{align*}
	which contradicts the anisotropy of $\varphi$.
	If we have $a\in -\D_F(2^m)=-\Gq_F(2^m)$, we have
	\begin{align*}
        		\Pfister{-1,\ldots, -1, a}&=\Pfister{-1,\ldots, -1}\perp -a\Pfister{-1,\ldots, -1}\\
        		&\cong\Pfister{-1,\ldots, -1}\perp \Pfister{-1,\ldots, -1}=\Pfister{\underbrace{-1,\ldots, -1}_{m+1}},
    	\end{align*}
	contradicting the maximality of $m$.
\end{proof}

By abuse of notation, we will henceforth often identify an element
$x\in F^\ast$ with its class in $F^\ast/ \D_F(2^m)$ for a given
$m\in\N$, and denote by $U_{n,m}$ the set of $(n-m)$-dimensional
subspaces of $F^\ast/ \D_F(2^m)$.

\begin{corollary}\label{UVRzuPfisterform}
	\begin{enumerate}[(a)]
		\item\label{cor:UVRzuPisterformParta} We have a well defined map 
	    	\begin{align*}
        			\varphi_{n,m}:U_{n,m}&\to \P_{n}(F)\\ 
			U&\mapsto\Pfister{-1, \ldots, -1, -x_1, \ldots, -x_{n-m}}
	    	\end{align*}
		where $\{x_1,\ldots, x_{n-m}\}$ is an arbitrary $\mathbb F_2$-basis of $U$.
		We then write 
		$$\pi_{U, m}=\varphi_{n,m}(U).$$
    		\item\label{cor:UVRzuPisterformPartb} Every
form in $\P_{n,m}(F)$ is in the image of $\varphi_{n,m}$.
		In particular, we have $|\P_{n,m}(F)|\leq |U_{n,m}|$.
	\end{enumerate}
\end{corollary}
\begin{proof}
	Let $x,y\in F^\ast$ be such that there is a $c\in\D_F(2^m)=\Gq_F(2^m)$
with $y=cx$.
	We then have
	\[\Pfister{-1,\ldots, -1,-y}=\Pfister{-1,\ldots, -1,-cx}
\cong\Pfister{-1,\ldots, -1,-x}.\]
	Furthermore,  any basis from a subspace $U$ can be transformed into any other basis of $U$ by successively replacing a subset $\{x,y\}$ by $\{y, x\}$ or by $\{x, xy\}$.
	Also, we have isometries $\Pfister{-x, -xy}\cong\Pfister{-x,-y}\cong\Pfister{-y,-x}$ for all $x,y\in F^\ast$.
	
	The map is thus well-defined and we have proved \eqref{cor:UVRzuPisterformParta}. 
    	Finally, for \eqref{cor:UVRzuPisterformPartb}, every form in $\P_{n,m}(F)$ has a preimage due to \Cref{LinearIndependentRep}.
\end{proof}

\begin{example} In \cite{MR0325656}, Cordes introduced
	so-called \emph{$\overline C$-fields}.  These are fields with
finite square class group, i.e., $q(F)<\infty$,  and such that $|\W(F)|=2^{q(F)}$.  
	Such fields will always have level $s(F)=1$ or $s(F)=2$.
	
	Over such fields, every subset of $F^\ast/F^{\ast2}$ is the value set of a unique anisotropic quadratic form over $F$ by \cite[Chapter XI. Theorem 7.19 (4)]{Lam2005}.
	In particular, we have $\D_F(\pi_{U,0})=U$ and $\D_F(\pi_{U, 1})=U\cup-U$.
	One now readily sees that each form in $\P_{n,m}(F)$ has exactly one preimage under $\varphi_{n,m}$.
	
	Note that $\overline C$-fields are examples of so-called \emph{rigid fields} whose \emph{Pfister number}, an invariant related to the symbol length, was studied by the second author in \cite{MR4554566}.
\end{example}

\begin{definition}\label{def:qmsmdm}
	Let $m$ be a nonnegative integer. For $m\geq 1$, we define
$$q_m=\dim_{\F_2}\left(D_F(2^m)/D_F(2^{m-1})\right).$$
If $q_m<\infty$, we have 
$$2^{q_m}=\left[\D_F\left(2^m\right):\D_F\left(2^{m-1}\right)\right].$$
For $m\geq 0$, we define
$$s_m=\dim_{\F_2}\left(D_F(\pm 2^m)/D_F(2^{m})\right)=
\log_2\left(\left[\D_F\left(\pm 2^m\right):\D_F\left(2^{m}\right)\right]\right)
\in\{ 0,1\}$$
and
$$d_m=\dim_{\F_2}\left(F^\ast/\pm \D_F\left(2^m\right)\right).$$
If $d_m<\infty$, we have
$$2^{d_m}=\left[F^*: \pm \D_F\left(2^{m}\right)\right].$$
\end{definition}

Note that $p(F)$ is finite if some $d_m$ is finite.
It is easy to calculate $s_m$ in terms of the level $s(F)$.

\begin{remark}\label{rem:calculateSm}
	\begin{enumerate}[(a)]
		\item\label{rem:calculateSm1} Let $F$ be a field of finite level $s(F)=2^s$ and let $\pi$ be the $m$-fold Pfister form $\Pfister{-1, \ldots, -1}$. 
    		Using the roundness of $\pi$, we see that for $m<s$, we have $\D_F(\pi)\cap-\D_F(\pi)=\emptyset$ and for $m\geq s$, we have $\D_F(\pi)=-\D_F(\pi)$. 
		In particular we have $s_m=1$ if $m<s$ and $s_m=0$ otherwise. 
    		\item If $F$ is a real field, we clearly have $s_m=1$ for all $m\in\N$.
\item  Note that $p(F)\leq 2^m$ iff $q_{m+1}=0$ iff $q_k=0$ for all $k\geq m+1$.
Thus, one readily sees that if $d_m<\infty$, then $p(F)\leq 2^{m+d_m+s_m}$.
	\end{enumerate}
\end{remark}

With this notation, we can formulate the following:

\begin{corollary}
Let $0\leq m\leq n$ be integers and let $F$ be a field with $d_m<\infty$.
	Then
   		 \[|\P_{n,m}(F)|\leq 2^{(n-m)(d_m+s_m-n+m)}\cdot \alpha_{n-m},\]
   	 where $\alpha_{n-m}$ is defined as in \Cref{cor:BoundNumberSubspaces}.
\end{corollary}
\begin{proof}
	This follows immediately from \Cref{UVRzuPfisterform} and \Cref{cor:BoundNumberSubspaces}.
\end{proof}

The rest of the section is now devoted to finding upper bounds for the $d_m$ in order to optimize the above bound for the number of Pfister forms.
We will see in the next result that we can express $d_m$ in terms of $s_m$ and the $q_k$ for $k\leq m$.

\begin{proposition}\label{DimFModPM2mQuad}
	Let $F$ be a field with $q(F)=2^d<\infty$. 
	We then have $$d_m=d-s_m-\sum\limits_{k=1}^mq_k$$
	for all $m\in\N$.
\end{proposition}
\begin{proof}
	For all $m\in\N$, we have 
	\begin{align*}
		2^d&=[F^\ast:F^{\ast2}]\\
		&=[F^\ast:\pm\D_F(2^m)]\cdot[\pm\D_F(2^m):\D_F(2^m)]\cdot\prod_{k=1}^m[\D_F(2^k):\D_F(2^{k-1})]\\
		&=2^{d_m}\cdot 2^{s_m}\cdot \prod_{k=1}^m2^{q_k}=2^{d_m+s_m+\sum\limits_{k=1}^mq_k}.
	\end{align*}
	The claim follows readily.
\end{proof}

We now want to collect estimates for the invariants $q_m$.
To do so we will use the following extension of a classical result from I.~Kaplansky.

\begin{theorem}\label{thm:Kaplansky}
	Let $F$ be a field with finite $p(F)$, i.e., $2^s\leq p(F)<2^{s+1}$ for some integer $s\geq 0$. For $k\in\{1,\ldots, s\}$, we have
		\[q_k\geq {s+1-k}.\]
	For $k=s+1$, we have $q_k\geq 1$ if $p(F)>2^s$ and $q_k=0$ if $p(F)=2^s$.
	For $k>s+1$, we have $q_k=0$. 
\end{theorem}
\begin{proof}
	For $k\in\{1,\ldots, s\}$ the proof can be found in \cite[Chapter XI. Kaplansky's Lemma 7.1]{Lam2005}.
	For $k=s+1$, we have $q_k=0$ if and only if $D_F(2^s)=D_F(2^{s+1})$ if and only if $p(F) = 2^s$.
	For $k>s+1$, we clearly have $D_F(2^k)=D_F(2^{k+1})$ and the assertion is clear.
\end{proof}

\begin{remark}
(i) If $F$ in the above theorem is nonreal, then $2^s\leq p(F)<2^{s+1}$ is equivalent
to $s(F)=2^s$, in which case $p(F)\in\{ 2^s,2^s+1\}$.

(ii) In the case $k=s+1$ and $p(F)>2^s$, $q_k$ can be arbitrary large. 
	To see this, consider the field $F_n=\F_3\laurent {t_1}\ldots\laurent {t_n}$.
	For all $n\geq 1$, we have $s(F_n)=2=2^1$ and $p(F_n)=3$, $\D_{F_n}(2^1)=\{1,2\}$ but $\D_{F_n}(2^2)=F_n^\ast$, so 
$q_2=n$.
\end{remark}

We are finally in a good position to give upper bounds for $d_m$ that only depend on $m$, on the Pythagoras number of $F$ and on whether $F$ is real or not.

\begin{corollary}\label{DimFSternModPM2mQuadrate}
Let $F$ be a field with $q(F)=2^d<\infty$ and let $s\geq 0$ be the integer with
$2^s\leq p(F)<2^{s+1}$.
	\begin{enumerate}[(a)]
		\item\label{bullet:DimFSternModPM2mQuadrate1} For $m< s$, we have
			\begin{align*}
				d_m\leq d-m\cdot\frac{2s-m+1}2-1.
			\end{align*}
		\item\label{bullet:DimFSternModPM2mQuadrate2} For $m = s$, we have
			\begin{align*}
				d_m\leq \begin{cases}
					d-s\cdot\frac{s+1}2,&\text{if }F\text{ is nonreal} \\
					d-s\cdot\frac{s+1}2-1,&\text{if }F\text{ is real}.
					\end{cases}
			\end{align*}
		\item\label{bullet:DimFSternModPM2mQuadrate3} For $m = s + 1$:
		\begin{enumerate}[(i)]
			\item if $F$ is nonreal, then $d_m=0$.
			\item\label{subbullet1} if $F$ is real with $p(F)=2^s$, we have
			 \begin{align*}
				d_m\leq d-s\cdot\frac{s+1}2-1.
			\end{align*}
			\item\label{subbullet2} if $F$ is real with $p(F)> 2^s$, we have 
				\begin{align*}
					d_m\leq d-s\cdot\frac{s+1}2-2.
				\end{align*}
		\end{enumerate}
		\item\label{bullet:DimFSternModPM2mQuadrate4} For $m > s+1$ we have
			 \begin{align*}
					d_m=0\text{ if }F\text{ is nonreal and }d_m\leq
						d-s\cdot\frac{s+1}2-1,&\text{ if }F\text{ is  real}.
			\end{align*}
	\end{enumerate}
\end{corollary}
\begin{proof}
	(\ref{bullet:DimFSternModPM2mQuadrate1}), (\ref{bullet:DimFSternModPM2mQuadrate2}): We use \Cref{DimFModPM2mQuad} and plug in the values obtained in \Cref{rem:calculateSm} and \Cref{thm:Kaplansky}. 
		In fact, we have
		\begin{align*}
			d_m&=d-s_m-\sum\limits_{k=1}^mq_k\\
			&\leq d-s_m-\sum\limits_{k=1}^m(s+1-k)\\
			&=d-s_m-m\cdot s-m+\frac{m(m+1)}2\\
			&=d-m\cdot\frac{2s-m+1}{2}-s_m
		\end{align*}
		and the result follows.\\
		The nonreal case of (\ref{bullet:DimFSternModPM2mQuadrate3}) follows since in this case, the $m$-fold Pfister form $\Pfister{-1,\ldots, -1}$ is clearly universal. 
		The other parts of (\ref{bullet:DimFSternModPM2mQuadrate3}) and (\ref{bullet:DimFSternModPM2mQuadrate4}) now follow using the same arguments as above. 
		The details are omitted and left to the reader
\end{proof}


As we already mentioned, $\pi\equiv a\pi\mod \I^{n+1}(F)$ for all $\pi\in \P_n(F), a\in F^\ast$, we immediately get the following result.

\begin{lemma}\label{BegrenzeDurchAnzahlPfister}
	For any field $F$, we have
	\[\sup_{\varphi\in \I^n(F)}\{\symlen_n(\varphi)\}\leq |\P_n(F)|.\]
\end{lemma}

Thus we could deduce upper bounds for the symbol length from our above calculations.
Nevertheless we do not work out the details since these bounds will obviously grow exponentially, while the upper bounds we will obtain in the next section grow polynomially.

\section{Bounding the Symbol Length}\label{Sec4}

\begin{proposition}\label{thm:ConstructBasis}
    There is a basis $\mathcal B$ of $F^\ast/F^{\ast2}$ with a filtration
    \[\mathcal B= \mathcal B_0\supseteq \mathcal B_1\supseteq \mathcal B_2\supseteq\ldots\]
    such that for each $m\in\N$, the classes of the elements in $\mathcal B_m$ form a basis of $F^\ast/\D_F(2^m)$.
\end{proposition}
\begin{proof}
    By basic linear algebra, we can find sets $\mathcal A_m$ for any $m\in\N$ such that $\mathcal A_m$ is a basis of $\D_F(2^m)/F^{\ast2}$ and such that we have $\mathcal A_m\subseteq \mathcal A_{m+1}$.
    We can then extend $\bigcup_{m\in\N}\mathcal A_m$ to a basis $\mathcal B$ of  the square class group $F^\ast/F^{\ast2}$.
    The result now follows since we can identify $\mathcal B_m:=\mathcal B\setminus\mathcal A_m$ with a basis of
	\begin{align*}
		(F^\ast/F^{\ast2})/(\D_F(2^m)/F^{\ast2})\cong F^\ast/ \D_F(2^m).
	\end{align*}
\end{proof}


\begin{remark}\label{dmsm}
    \begin{enumerate}[(a)]
    	\item The union $\bigcup_{m\in\N}\mathcal A_m$ in the proof of \Cref{thm:ConstructBasis} is a basis itself of $\D_F(\infty)/F^{*2}$
and thus equal to $\mathcal B$ if and only if every element is a sum of squares, i.e.,
if and only if $F$ is nonreal.
    	\item With our notation from \Cref{def:qmsmdm}, we have $|\mathcal B_m|=d_m+s_m$.
	\end{enumerate}
\end{remark}

The following easy lemma will be used frequently in the sequel.

\begin{lemma}\label{lem:ReplaceBySumOf2}
	For all $x,y,  x_2,\ldots, x_k\in F^\ast$, we have 
	\begin{align*}
		\Pfister{xy,x_2,\ldots, x_k}\equiv\Pfister{x,x_2,\ldots, x_k}+\Pfister{y,x_2,\ldots, x_k}\mod \I^{k+1}(F).
	\end{align*}
\end{lemma}
\begin{proof}
	From the Witt equivalence $\Pfister{xy}=\Pfister x + x\Pfister y$, we obtain 
	\begin{align*}
		\Pfister{xy,x_2,\ldots, x_k}&=\Pfister{x, x_2, \ldots, x_k} + x \Pfister{y, x_2,\ldots, x_k}\nonumber\\
		&\equiv\Pfister{x,x_2,\ldots, x_k}+\Pfister{y,x_2,\ldots, x_k}\mod \I^{k+1}(F)
	\end{align*}
\end{proof}

\begin{proposition}\label{GuteDarstellungVonFormen}
	Let $\mathcal B$ be a basis of $F^\ast/F^{\ast2}$ as in \Cref{thm:ConstructBasis}. 
    For $\varphi\in \I^n(F)$, there are finite subsets $C_m\subseteq \mathcal P_m(\mathcal B_{n-m})$ for $m\in\{0,\ldots, n\}$ such that we have
	\begin{align*}
		\varphi\equiv\sum\limits_{m=0}^n\sum\limits_{U\in C_m}\pi_{\spn(U), n-m}\mod \I^{n+1}(F).
	\end{align*}
\end{proposition}
\begin{proof}
	We consider a representation

	\[\varphi\equiv \pi_1+\ldots+\pi_k\mod \I^{n+1}(F)\]
	such that each $\pi_j$ is written with as many slots equal to $-1$ as
possible, i.e. if $\pi_j \in \P_{n,m}(F)$ for some $m\in\N$, then exactly $m$ slots
of $\pi_j$ are equal to $-1$ and the other $n-m$ slots can be chosen to be of
shape $-\prod_{i=1}^\ell b_i$ with $b_i\in \mathcal B_m$ and
suitable $\ell$, see \Cref{UVRzuPfisterform}.

	If we have $\pi_j=\Pfister{x_1, \ldots, x_{n-1}, x_n} $,
where $x_n=-b_1\cdot\ldots\cdot b_\ell=(-1)^{\ell-1}\prod_{i=1}^\ell(-b_i)$ with $b_1,\ldots, b_\ell\in\mathcal B_m$
and $\ell\geq 2$, for some $j\in\{1,\ldots, k\}$,  then by \Cref{lem:ReplaceBySumOf2},
we can replace $\pi_j$ by
the sum
$$\sum\limits_{i=1}^\ell \Pfister{x_1, \ldots, x_{n-1}, -b_{i}}
+\Pfister{x_1, \ldots, x_{n-1}, (-1)^{\ell-1}}.$$

	If we have $\pi_j\in \P_{n,m}(F)$ for some $m\in\N$, the resulting Pfister 
forms are hyperbolic or will lie in $\P_{n,m'}(F)$ for some $m'\geq m$.
	We can clearly omit the hyperbolic Pfister forms in our representation and 
thus now consider any of the resulting Pfister forms lying in $\P_{n,m'}(F)$
for some $m'\geq m$
	
    	If we have $m'=m$, we can repeat the above until the considered slot in any 
newly introduced Pfister form lies in $\mathcal B_m$. 
	If we have $m'>m$, we first substitute the representation of this form with 
an isometric Pfister form with $m'$ slots equal to $-1$ and all other slots being 
products of $-1$ and of elements lying in $\mathcal B_{m'}$.
	We can then apply the above procedure to this new form. 
	
	By repeating this process for all of the $\pi_j, j\in\{1,\ldots, k\}$ and all 
slots that neither lie in the appropriate $\mathcal B_m$ (up to a factor $-1$)
nor are equal to $-1$, and possibly rearranging the terms,
we will eventually obtain a representation as desired.
\end{proof}

\begin{remark}\label{AngepassteSummenGrenze}
	If $F$ is nonreal of level $s(F)=2^s$, it is obviously enough to restrict the outer sum occuring in \Cref{GuteDarstellungVonFormen} to $m\geq \max\{ 0,n-s\}$.
\end{remark}

We can directly deduce the following upper bound for the symbol length from \Cref{GuteDarstellungVonFormen}.

\begin{corollary}\label{SchrankeBasisStandard}
	Let $F$ be a field with $q(F)=2^d<\infty$,  and let $n\in\N$ be an integer.
	\begin{enumerate}[(a)]
		\item If $F$ is nonreal with $s(F)=2^s$, we have
			\begin{align*}
				\symlen_n(F)\leq \sum\limits_{m=0}^{\min\{s,n\}} \binom {d_m+s_m}{n-m}
			\end{align*}
		\item If $F$ is  real, we have 
			\begin{align*}
				\symlen_n(F)\leq \sum\limits_{m=0}^{n} \binom {d_m+s_m}{n-m}
			\end{align*}
	\end{enumerate}
\end{corollary}

We will now present two strategies for adding certain Pfister forms that appear in \Cref{GuteDarstellungVonFormen} in order to obtain better bounds for the symbol length.

\begin{corollary}\label{SchrankeBasisLink}
	Let $F$ be a field with $q(F)=2^d<\infty$, and let $n\in\N$.
	\begin{enumerate}[(a)]
		\item If $F$ is nonreal with $s(F)=2^s$, we have
			\begin{align*}
				\symlen_n(F)\leq\sum\limits_{m=0}^{\left\lfloor \frac{\min\{s,n\}}2\right\rfloor} \binom {d_{2m}+s_{2m}}{n-{2m}}.
			\end{align*}
		\item If $F$ is  real, we have 
			\begin{align*}
				\symlen_n(F)\leq \sum\limits_{m=0}^{\left\lfloor \frac n2\right\rfloor} \binom {d_{2m}+s_{2m}}{n-2m}.
			\end{align*}
	\end{enumerate}
\end{corollary}
\begin{proof}
	We write $\varphi$ as in \Cref{GuteDarstellungVonFormen}
and use the notation of the statement of this proposition.
	Fix a subset $A=\{x_1,\ldots, x_{n-2m}\}\in \Pcal_{n-2m}(\mathcal B_{2m})$
for some $m\geq1$.
	We now consider the set
$C\subseteq C_{n-2m+1}\subseteq \mathcal P_{n-2m+1}(\mathcal B_{2m-1})$
consisting of all $U\in C_{n-2m+1}$ with $A\subseteq U$.
	Let $a_1,\ldots, a_k\in \mathcal B_{2m-1}$ be pairwise different such that 
		\[C=\{A\cup\{a_j\}\mid j=1,\ldots, k\}.\]
	The Pfister forms
	\begin{align*}
		\Pfister{-1}^{\otimes 2m-1}\otimes\Pfister{-x_1, \ldots, -x_{n-2m}}
\otimes\Pfister{-a_j}, \text{ for }j\in\{1,\ldots, k\}
	\end{align*}
	then all occur in the given representation of $\varphi\mod\I^{n+1}(F)$.
	It may further happen that the Pfister form
	\begin{align*}
		\Pfister{-1}^{\otimes 2m}\otimes\Pfister{-x_1, \ldots, -x_{n-2m}}
	\end{align*}
	(corresponding to the set $A$ from the beginning of the proof) occurs
in the representation. 
	By \Cref{lem:ReplaceBySumOf2} the sum (modulo $\I^{n+1}(F)$)
of these Pfister forms is
equivalent to 
		\[\Pfister{-1}^{\otimes 2m-1}\otimes\Pfister{-x_1, \ldots, -x_{n-2m}}
\otimes \Pfister{\epsilon a_1\cdot \ldots \cdot a_k}\]
where $\epsilon\in\{\pm 1\}$ depends on the parity of $k$ and on	
whether $\Pfister{-1}^{\otimes 2m}\otimes\Pfister{x_1, \ldots, x_{n-2m}}$
has to be included or not.
	We will sum up all Pfister forms in the given representation
except those in $\P_{n,0}(F)$ if we apply this procedure for all possible
$A$ for all $m\geq1$. 
	There are at most $\binom{d_0+s_0}{n}$ elements in $\P_{n,0}(F)$,
and for $m\geq 1$, there are at most $\binom{d_{2m}+s_{2m}}{n-2m}$ possible
choices of a subset $A$ as above. 
	Thus the claim follows.
\end{proof}

\begin{theorem}\label{PolynomialBound}
	Let $F$ be a field with $q(F)=2^d<\infty$ and let $n\in\N$ be an integer. 
	\begin{enumerate}[(a)]
		\item\label{PolynomialBound1} If $F$ is nonreal with $s(F)=2^s$, we have
			\begin{align*}
				\symlen_n(F)\leq \sum\limits_{m=0}^{\min\{s,n-1\}} \sum\limits_{r=0}^{\left\lfloor\frac{n-m-1}2\right\rfloor} \binom{\left\lfloor \frac{d_m+s_m}2\right\rfloor}{2r}\binom{\left\lfloor \frac{d_m+s_m+1}2\right\rfloor}{n-m-1-2r}
			\end{align*}
		\item If $F$ is  real, we have 
			\begin{align*}
				\symlen_n(F)\leq \sum\limits_{m=0}^{n-1} \sum\limits_{r=0}^{\left\lfloor \frac{n-m-1}2\right\rfloor} \binom{\left\lfloor \frac{d_m+s_m}2\right\rfloor}{2r}\binom{\left\lfloor \frac{d_m+s_m+1}2\right\rfloor}{n-m-1-2r}
			\end{align*}
	\end{enumerate}
\end{theorem}
\begin{proof}
    	For all $m\in\{0,\ldots, n\}$, we decompose $\mathcal B_m$ into $\mathcal B_m=\mathcal A_m\cup\mathcal A_m'$ with $\mathcal A_m\cap\mathcal A_m'=\emptyset$ and $|\mathcal A_m|=\left\lfloor\frac{d_m+s_m}{2}\right\rfloor$ for some suitably chosen subsets $\mathcal A_m \subseteq \mathcal B_m$. 
    	We then have (see also \Cref{dmsm}(b))
		\[|A_m'|=d_m+s_m-\left\lfloor\frac{d_m+s_m}2\right\rfloor=\left\lfloor\frac{d_m+s_m+1}{2}\right\rfloor.\]
    	We further consider
    		\[\Pcal_{m,r}=\{X\subseteq \mathcal A_m\mid |X| = 2r\} ~~~\text{ and }~~~ \Pcal_{m,r}'=\{Y\subseteq \mathcal A_m'\mid |Y|=n-m-1-2r\}.\]
    	We clearly have
    		\[|\Pcal_{m,r}|=\binom{\left\lfloor\frac{d_m+s_m}{2}\right\rfloor}{2r} ~~~\text{ and }~~~ |\Pcal_{m,r}'|=\binom{\left\lfloor\frac{d_m+s_m+1}{2}\right\rfloor}{n-m-1-2r}.\]
    	Let now $\varphi\in \I^n(F)$ and consider a representation 
	\begin{align}\label{eq:MainEquivalence}
    		\varphi\equiv\sum\limits_{m=0}^n\sum\limits_{U\in C_m}\pi_{U, n-m}\mod \I^{n+1}(F)
	\end{align}
    	as in \Cref{GuteDarstellungVonFormen}. 
	Let $X\in \Pcal_{m,r}$ and $Y\in \Pcal_{m,r}'$. 
	Then all forms $\pi\in \P_{n, m}(F)\cup \P_{n, m+1}(F)$ that occur
in this representation and that have the elements of $X\cup Y$ as slots
(up to a factor $-1$)
can be replaced by a single Pfister form using \Cref{lem:ReplaceBySumOf2}
just as in the proof of \Cref{SchrankeBasisLink}.
	
	Conversely, let some $\pi=\pi_{U,m}$ from this representation be given
and let $B=\{-b_1,\ldots, -b_{n-m}\}$ be the set of slots $\neq-1$ of $\pi$.
	Let $r$ be maximal such that there is an $X\in \Pcal_{m,r}$ with
$X\subseteq -B$.
	We clearly have 
		\[|(-B\setminus X)\cap \mathcal A_m'|\in\{n-m-2r-1, n-m-2r\}.\]
	There is thus a subset $Y\subseteq \mathcal A_m'$ with cardinality
$|Y|=n-m-2r-1$ such that $X\cup Y\subseteq -B$.
	
	Hence, for all Pfister forms of type $\pi_{U, m}$ for a suitable $U$, 
there is an $r\in\N$, $X\in \Pcal_{m,r}$ and $Y\in \Pcal_{m,r}'$ such that
all slots except one are either $-1$ or lie (up to a factor $-1$) in $X\cup Y$.
	Therefore, our above strategy to replace certain Pfister forms by a
single one can be applied to all forms on the right hand side
in \eqref{eq:MainEquivalence}.
	
	The symbol length $\symlen_n(\varphi)$ can thus be bounded by
the cardinality of 
	\[\{X\times Y\mid \exists m\in\N, r\in\N: X\in \Pcal_{m,r}, Y\in \Pcal_{m,r}'\}\]
	which is given by the respective sums of the theorem.
\end{proof}

\begin{lemma}\label{BinKoeffAlsPoly}
	Let $k,m\in\N$ be nonnegative integers with $m\geq2$. 
    Then 
	\[\sum\limits_{r=0}^{\left\lfloor\frac{n-m-1}2\right\rfloor}\binom{k}{2r}\binom{k}{n-m-1-2r}\text{ and }\sum\limits_{r=0}^{\left\lfloor\frac{n-m-1}2\right\rfloor}\binom{k}{2r}\binom{k+1}{n-m-1-2r}\]
	are polynomials in $k$ of degree $n-m-1$ with leading coefficient $\frac{2^{n-m-1}}{2(n-m-1)!}$.
\end{lemma}
\begin{proof}
	For all $r$, the expression $\binom{k}{2r}$ is zero or a polynomial in $k$ of degree $2r$ with leading coefficient $\frac{1}{(2r)!}$.
Similarly $\binom{k}{n-m-1-2r}$ and $\binom{k+1}{n-m-1-2r}$ are zero or polynomials in $k$ of degree $n-m-1-2r$ with leading coefficient $\frac1{(n-m-1-2r)!}$. 
    Thus, in both cases, the products are zero or polynomials in $k$ of degree
    \[2r+n-m-1-2r=n-m-1\] 
    with leading coefficient
    \[\frac1{(2r)!(n-m-1-2r)!}.\] 
    The sum then is a polynomial of degree $n-m-1$ with leading coefficient
	\begin{align*}
		\sum\limits_{r=0}^{\left\lfloor\frac{n-m-1}2\right\rfloor}\frac1{(2r)!(n-m-1-2r)!}&=\frac{1}{(n-m-1)!}\sum\limits_{r=0}^{\left\lfloor\frac{n-m-1}2\right\rfloor}\frac{(n-m-1)!}{(2r)!(n-m-1-2r)!}\\
		&=\frac{1}{(n-m-1)!}\sum\limits_{r=0}^{\left\lfloor\frac{n-m-1}2\right\rfloor}\binom{n-m-1}{2r}\\
		&=\frac1{(n-m-1)!}2^{n-m-2}\\
		&=\frac{2^{n-m-1}}{2(n-m-1)!}.
	\end{align*}
\end{proof}

\begin{corollary}\label{kahn-bound}
	Let $F$ be a field with $q(F)=2^d<\infty$ and let $n\geq 1$ be an integer.
    Then $\symlen_n(F)$ is bounded from above by 
	\[\frac{d^{n-1}}{2(n-1)!}+f_n(d)\]
	where $f_n$ is a polynomial of degree at most $n-2$.
\end{corollary}
\begin{proof}
	This follows from the above by plugging in $\left\lfloor\frac{d_0+s_0}2\right\rfloor$ for $k$ in \Cref{BinKoeffAlsPoly} and noting that 
the sums in \Cref{PolynomialBound} can be bounded from above
by replacing any $d_m+s_m$ by $d_0+s_0=d$.
\end{proof}

\begin{remark}
	In the above corollary, we rediscovered \cite[Proposition 2.3h)]{Kahn}
for fields with finite square class number $2^d$ in a way that allows us to
determine the polynomial explicitly. 
    	Of course, sticking with  $d_m+s_m$ in \Cref{PolynomialBound}
rather than replacing them by $d$
will yield better bounds in general but also makes their computation more difficult.
\end{remark}

\begin{example}\label{example-n=3}
We now close the article with a comparison of these bounds for level $s(F)=2^s$ with $s\in\{0,1,2\}$. (Recall that conjecturally, if $s(F)<\infty$ then $q(F)<\infty$ implies $s(F)\leq 4$.)
Of interest in this context are certain fields of iterated Laurent series that realize a
given value $q(F)$ for the respective value $s(F)$, such as
\begin{itemize}
\item $s=1$, $q=2^d$ ($d\geq 0$): $\CC(\!(t_1)\!)\ldots (\!(t_d)\!)$;
\item $s=2$, $q=2^d$ ($d\geq 1$): $\F_3(\!(t_1)\!)\ldots (\!(t_{d-1})\!)$;
\item $s=4$, $q=2^d$ ($d\geq 3$): $\Q_2(\!(t_1)\!)\ldots (\!(t_{d-3})\!)$.
\end{itemize}
(Note that for any field $F$ with $s(F)=4$, one has $q(F)\geq 8$, see \cite[Theorem 2.7]{MR0323716}.)

    Before giving some explicit estimates, we want to recall the upper bounds that we found above.

	First of all, we want to collect upper bounds for $d_m$ for some cases. 
 	In the following table, we have summarized the upper bounds we know for $d_m$ obtained in \Cref{DimFSternModPM2mQuadrate}:
	
	\begin{tabular}{|c||c|c|c|}
		\hline
		\diagbox{$m$}{$s$} & 0 & 1 & 2\\
		\hline
		\hline
		0 & $d$ & $d-1$ & $d-1$\\
		\hline
		1 & $0$ & $d-1$ & $d-3$\\
		\hline
		2 & $0$ & $0$ & $d-3$\\
		\hline
	\end{tabular}
	
	For $s_m$ we have the following values according to \Cref{rem:calculateSm} (\ref{rem:calculateSm1}):
	
	\begin{tabular}{|c||c|c|c|}
		\hline
		\diagbox{$m$}{$s$} & 0 & 1 & 2\\
		\hline
		\hline
		0 & $0$ & $1$ & $1$\\
		\hline
		1 & $0$ & $0$ & $1$\\
		\hline
		2 & $0$ & $0$ & $0$\\
		\hline
	\end{tabular}

	Now we plug in the upper bounds for $d_m$ and $s_m$ in all the bounds for the symbol length that we have found so far. 
    It is obvious that this will also yield upper bounds for the symbol length.
	
	\begin{itemize}
		\item \Cref{SchrankeBasisStandard} yields
			$$\symlen_n(F)\leq\begin{cases}
				\binom dn, &\text{ if }s=0\\
				\binom{d}n+\binom{d-1}{n-1}, &\text{ if }s=1\\
				\binom{d}n+\binom{d-2}{n-1}+\binom{d-3}{n-2}, &\text{ if }s=2.
			\end{cases}$$
		\item \Cref{SchrankeBasisLink} yields
			$$\symlen_n(F)\leq\begin{cases}
				\binom d{n}, &\text{ if }s=0\\
				\binom{d}{n}, &\text{ if }s=1\\
				\binom{d}{n}+\binom{d-3}{n-2}, &\text{ if }s=2.
			\end{cases}$$
		\item Since we do not have a closed formula for the upper bound obtained in \Cref{PolynomialBound}, we will only write down the case $n=3$ explicitly. 
        We obtain the following upper bounds:\\
			\begin{tabular}{|c||c|c|c|}
				\hline
				$s$ & 0 & 1 & 2\\
				\hline
				\hline
				bound & 
				$\binom{\left\lfloor\frac{d+1}{2}\right\rfloor}{2}+ \binom{{\left\lfloor\frac{d}{2}\right\rfloor}}{2}$ & 
				$\binom{\left\lfloor\frac{d+1}{2}\right\rfloor}{2}+\binom{\left\lfloor\frac{d}{2}\right\rfloor}{2}+\left\lfloor\frac{d}{2}\right\rfloor$ & 
				$\binom{\left\lfloor\frac{d+1}{2}\right\rfloor}{2}+ \binom{\left\lfloor\frac{d}{2}\right\rfloor}{2}+{\left\lfloor\frac{d-1}{2}\right\rfloor}+1$ \\
if $d=2k$ & $k^2-k$ & $k^2$ & $k^2$\\
if $d=2k+1$ & $k^2$ & $k^2+k$ & $k^2+k+1$\\				
\hline
			\end{tabular}
	\end{itemize}
Let us illustrate these upper bounds for the $3$-symbol lengths for certain values
of $d$:
	\begin{itemize}
		\item For $s=0$, i.e. $s(F)=1$:\\
		\begin{tabular}{|c|c|c|c|c|c|c|c|c|}
			\hline
			& $d=4$ & $d=5$ & $d=7$ & $d=10$\\
			\hline
			\Cref{SchrankeBasisStandard} & 4 & 10 & 35 & 120\\
			\hline
			 \Cref{SchrankeBasisLink}  & 4 & 10 & 35 & 120\\
			\hline
			\Cref{PolynomialBound} & 2 & 4 & 9 & 20\\
			\hline
		\end{tabular}~\\
		\item For $s=1$, i.e. $s(F)=2$:\\
		\begin{tabular}{|c|c|c|c|c|c|c|c|c|}
			\hline
			& $d=4$ & $d=5$ & $d=7$ & $d=10$  \\
			\hline
			\Cref{SchrankeBasisStandard} & 7 & 16 & 50 & 156 \\
			\hline
			 \Cref{SchrankeBasisLink} & 4 & 10 & 35 & 120\\
			\hline
			\Cref{PolynomialBound} & 4 & 6 & 12 & 25 \\
			\hline
		\end{tabular}~\\
		\item For $s=2$, i.e. $s(F)=4$:\\
		\begin{tabular}{|c|c|c|c|c|c|c|c|c|}
			\hline
			& $d=4$ & $d=5$ & $d=7$ & $d=10$ \\
			\hline
			\Cref{SchrankeBasisStandard} & 6 & 15 & 49 & 155\\
			\hline
			 \Cref{SchrankeBasisLink} & 5 & 12 & 39 & 127\\
			\hline
			\Cref{PolynomialBound} & 4 & 7 & 13 & 25\\
			\hline
		\end{tabular}~\\
	\end{itemize}
	
	It comes as no suprise that the bound obtained in \Cref{PolynomialBound} yields the best values.
	
	We now return to the case in which $F$ is an iterated Laurent extension of $\CC, \mathbb F_3, \Q_2$ with $d=4$ as above. 
It is straightforward to see that $F$ is 3-linked in these cases, i.e. for each pair of 3-fold Pfister forms $\pi_1, \pi_2$ over $F$, there are $a, b, c_1, c_2\in F^\ast$ with $\pi_1\cong \Pfister{a, b, c_1}, \pi_2\cong \Pfister{a, b, c_2}$. 
In particular, we have $\symlen_3(F)=1$, strictly subceeding the upper bounds above.
\end{example}

\begin{example}
	Let us consider the special case of a field $F$ with level $s(F)=1=2^0$.
	In \Cref{PolynomialBound} (\ref{PolynomialBound1}), we then have $m=0$ and $d_0+s_0=d$, so the upper bound simplifies to
	\begin{align}\label{eq:boundForLevel0}
		\symlen_n(F)\leq \sum\limits_{r=0}^{\left\lfloor\frac{n-1}2\right\rfloor} \binom{\left\lfloor \frac{d}2\right\rfloor}{2r}\binom{\left\lfloor \frac{d+1}2\right\rfloor}{n-1-2r}.
	\end{align}
	We now have a closer look at the case $n>d$.
	By \cite[Chapter XI. Kneser's Lemma 6.5]{Lam2005} there are no anisotropic forms of dimension $>2^d$.
	In particular, we have $\I^n(F)=0$ and thus $\symlen_n(F)=0$.
	For $n\geq d+2$, one readily sees that the right hand side of \eqref{eq:boundForLevel0} also equals 0.
	For $n=d+1$, the right hand side of \eqref{eq:boundForLevel0} equals $1$ if $d\equiv0,1\mod4$ and equals $0$ if $d\equiv2,3\mod4$.
	The details of these straightforward calculations are left to the reader.
\end{example}


\bibliographystyle{alpha}
\bibliography{literatur}


\bibliographystyle{alpha}
\bibliography{literatur}

\end{document}